\documentclass[12pt, twoside, leqno]{amsart}
\usepackage{amsmath}
\usepackage{amsfonts}
\usepackage{amssymb}
\usepackage{amsthm}
\usepackage{url}
\usepackage{latexsym}
\usepackage{bbm}
\usepackage{enumerate}
\usepackage{booktabs}
\usepackage[T1]{fontenc}

\newcommand{\R}{\mathbb{R}}
\newcommand{\Z}{\mathbb{Z}}
\newcommand{\Rn}{\mathbb{R}^n}

\newcommand{\Ar}{\Arrowvert}
\newcommand{\D}{\, \mathrm{d}}

\newcommand{\comment}[1]{}

\newcommand{\intav}{-\hskip -13pt\int}

\DeclareMathOperator{\bmo}{BMO}

\def\mvint_#1{\mathchoice%
          {\mathop{\kern 0.2em\vrule width 0.6em height 0.69678ex depth -0.58065ex
                  \kern -0.8em \intop}\nolimits_{\kern -0.4em#1}}%
          {\mathop{\kern 0.1em\vrule width 0.5em height 0.69678ex depth -0.60387ex
                  \kern -0.6em \intop}\nolimits_{#1}}%
          {\mathop{\kern 0.1em\vrule width 0.5em height 0.69678ex depth -0.60387ex
                  \kern -0.6em \intop}\nolimits_{#1}}%
          {\mathop{\kern 0.1em\vrule width 0.5em height 0.69678ex depth -0.60387ex
                  \kern -0.6em \intop}\nolimits_{#1}}}

\def\sqr#1#2{{\vcenter{\hrule height.#2pt\hbox{\vrule
    width.#2pt height#1pt\kern#1pt\vrule width.#2pt}\hrule height.#2pt}}}
\def\square{\mathchoice\sqr{5.5}4 \sqr{5.0}4 \sqr{4.8}3 \sqr{4.8}3}
\def\XXint#1#2#3{{\setbox0=\hbox{$#1{#2#3}{\int}$}
     \vcenter{\hbox{$#2#3$}}\kern-.5\wd0}}

\theoremstyle{plain}

\numberwithin{equation}{section}

\theoremstyle{definition}

\theoremstyle{remark}

\title{Parabolic John-Nirenberg spaces}
\author{Lauri Berkovits}
\date{\today}

\begin{document}

\keywords{John-Nirenberg lemma, parabolic BMO}
\subjclass[2000]{46E30}

\begin{abstract}
We introduce a parabolic version of John-Nirenberg space with exponent $p$ and 
show that it is contained in local weak-$L^p$ spaces.
\end{abstract}

\maketitle

\section{Introduction}

In the classical paper of F. John and L. Nirenberg \cite{JoNi}, where 
functions of bounded mean oscillation ($\bmo$) were introduced, 
they also studied a class satisfying a weaker $\bmo$ type condition
\begin{equation*}
K_f^p:=\sup_{\{Q_j\}_j}\sum_j |Q_j| \left( \mvint_{Q_j} |f-f_{Q_j}| \D x \right)^p <\infty,
\end{equation*}
where the supremum is taken over all partitions $\{Q_j\}_j$ 
of a given cube $Q_0$ into pairwise non-overlapping subcubes.
The functional $f \mapsto K_f$ defines a seminorm and the class of functions satisfying 
$K_f<\infty$, which we denote by $JN_p(Q_0)$ for John-Nirenberg space with exponent $p$,
 can be seen as a generalization 
of $\bmo$. Indeed, $\bmo$ is obtained as the limit case of $JN_p$ in the sense that
$$
\lim_{p \to \infty} K_f=\sup_{Q \subseteq Q_0} \mvint_Q |f-f_Q| \D x = ||f||_{\bmo(Q_0)}.
$$
In contrast to the exponential integrability of $\bmo$ functions, 
$K_f<\infty$ implies that $f$ belongs to the space weak-$L^p(Q_0)$.
This was already observed by John and Nirenberg. Precisely, they showed that for 
$\lambda >0$, we have
\begin{equation*}
|\{x \in Q_0: |f(x)-f_{Q_0}| > \lambda \}| \leq C\left(\frac{K_f}{\lambda}\right)^p,
\end{equation*}
where the constant $C$ depends on $n$ and $p$. 
Simpler proofs and generalizations 
have appeared in \cite{AaBeKaYu, BaJiMa, FrPeWh,GiMa, Gi,Heik1, JiMa, JoNi, LePe, MaPe,To}. In this note we show that an analogous result 
holds in the context of parabolic $\bmo$ spaces.

\section{Parabolic John-Nirenberg space}

We shall introduce some notation and terminology. 
Given an Euclidean cube $Q=\prod_{i=1}^n [a_i,a_i+h]$, 
we define the forward in time translation
$$
Q^{+}:=\prod_{i=1}^{n-1} [a_i,a_i+h] \times [a_n+h,a_n+2h].
$$ 
Moreover, we use the notation $Q^{+,2}:=(Q^+)^+$. 
We write $f \in \bmo^+(\Rn)$, if we have
\begin{equation}\label{bmoplus}
\Ar f \Ar_{\bmo^+(\Rn)}:= \sup_{Q} \mvint_Q (f-f_{Q^+})^+\D x < \infty,
\end{equation}
where the supremum is taken over all cubes in $\Rn$ 
with sides parallel to the coordinate axes. It should be observed that 
despite the notation, the quantity defined by \eqref{bmoplus} is not actually a norm.

The one-dimensional $\bmo^+(\R)$ class was first introduced by F. J. Mart\'in-Reyes 
and A. de la Torre \cite{MaTo1}, who showed that this class possesses many properties similar to 
the standard $\bmo$ space. Even though steps towards a multidimensional theory has been taken (see \cite{Be}), a 
satisfactory theory has only been developed in dimension one. 
In the classical elliptic setting, one of the cornerstones of theory of $\bmo$ functions is the 
celebrated John-Nirenberg inequality, which shows that logarithmic growth is the maximum possible 
for a $\bmo$ function. 
A corresponding result holds for the class $\bmo^+(\R)$, and a slightly weaker version 
of this result for $\bmo^+(\Rn)$ was obtained in \cite{Be}.

In this setting we define John-Nirenberg spaces as follows. 
We write $f \in$ $JN_p^+(\Rn)$ if 
\begin{equation}\label{jnpplus}
(K_f^+)^p:=\sup_{\{Q_j\}_j}
\sum_j |Q_j| \left( \intav_{Q_j \cup Q_j^+} (f-f_{Q_j^{+,2}})^+ \D x \right)^p <\infty,
\end{equation}
where the supremum is taken over countable families $\{Q_j\}$ of pairwise non-overlapping cubes 
satisfying $\sum_j |Q_j| < \infty$. 
The definition is reasonable in the sense that the $\bmo^+(\Rn)$ condition may be seen as the limit 
case of \eqref{jnpplus} as $p \to \infty$. Precisely, 
$$
\lim_{p \to \infty} K_f^+ =\sup_Q \intav_{Q \cup Q^+} (f-f_{Q^{+,2}})^+ \D x,
$$
where the quantity on the right-hand side is equivalent (up to a multiplication by a universal constant) to the $\bmo^+$ norm of $f$, 
defined by \eqref{bmoplus}.

The following theorem is a parabolic version of the weak distribution inequality of John 
and Nirenberg.

{\bf Theorem.}
\emph{ Assume $f \in JN_p^+(\Rn)$. Then, for every cube $Q_0$ and $\lambda>0$, we have
\begin{equation}\label{theorem}
|\{x \in Q_0: (f(x)-f_{Q_0^{+,2}})^+ > \lambda \}| \leq C\left(\frac{K_f^+}{\lambda}\right)^p,
\end{equation}
where $C$ only depends on $n$ and $p$.}

\section{Proof of the theorem}

We follow the argument used in \cite{AaBeKaYu}. 
Given a non-negative $f$ and a cube $Q_0$,
denote by $\Delta=\Delta(Q_0)$ the family of all dyadic subcubes
obtained from $Q_0$ by repeatedly bisecting the sides into 
two parts of equal length.
We shall make use of the ``forward in time dyadic maximal function'' defined by
$$
M_{Q_0}^{+,d}f(x):=\sup_{\substack{Q \in \Delta \\ x \in Q}} \intav_{Q^+} f \D x.
$$
A standard stopping-time argument shows that we have 
$$
\{x \in Q_0: M_{Q_0}^{+,d}f(x)>\lambda\}=\bigcup_j Q_j,
$$
where $Q_j$'s are the maximal dyadic subcubes of $Q_0$ satisfying
\begin{equation}\label{p1}
\intav_{Q_j^+} f \D x >\lambda. 
\end{equation}
Maximality implies that the cubes $Q_j$ are pairwise non-overlapping.
Moreover, if $\lambda \geq f_{Q_0^+} $, then $Q_0$ doesn't satisfy \eqref{p1}. 
Consequently, in this case every $Q_j$ is contained in a 
larger dyadic subcube $Q_{j^-}$ of $Q_0$ which does not satisfy 
\eqref{p1}. Since $Q_j^{+,2} \subset Q_{j^-}^+$, we conclude
\begin{equation}\label{p2}
\intav_{Q_j^{+,2}} f \D x \leq 2^n\lambda,
\end{equation}
provided $\lambda \geq f_{Q_0^+}$.
Standard arguments imply a weak type estimate for $M_{Q_0}^{+,d}$. Indeed, we have
$$
|\{x \in Q_0: M_{Q_0}^{+,d}f(x)>\lambda\}|=\sum_j |Q_j|.
$$
While the cubes $Q_j$ are non-overlapping, the cubes $Q_j^+$ may not be.
Let us replace $\{Q_j^+\}_j$ by the maximal non-overlapping subfamily 
$\{\widetilde{Q}_j^+\}_j$ which we form by collecting those $Q_j^+$ which 
are not properly contained in any other $Q_{j'}^+$.
Maximality of $\{\widetilde{Q}_{j}^+\}_{j}$ enables us to partition the family 
$\{Q_j\}_j$ as follows. Given $\widetilde{Q}_{j}^+$, we 
define $I_{j}:=\{ i : Q_{i}^+ \subseteq \widetilde{Q}_{j}^+\}$, and we may write 
$\{Q_j\}_j=\bigcup_j \{Q_i: i \in I_j\}.$
Now, whenever $i \in I_{j},$ we have $Q_{i}\subseteq \widetilde{Q}_j \cup\widetilde{Q}_j^+$ 
and we get the estimate
\begin{align*}
\sum_{j}|Q_{j}|&=\sum_{j} \sum_{i \in I_{j}}|Q_{i}| \\ &
\leq 2\sum_{j} |\widetilde{Q}_{j}^+| \\ &
\leq \frac{2}{\lambda} \int_{Q_0 \cup Q_0^+} f \D x.
\end{align*}
Combining the previous estimates, we arrive at
\begin{equation}\label{p3}
|\{x \in Q_0: M_{Q_0}^{+,d}f(x)>\lambda\}| \leq \frac{2}{\lambda} \int_{Q_0 \cup Q_0^+} f \D x.
\end{equation}

We begin by proving the following good $\lambda$ inequality for 
the forward in time dyadic maximal operator.

{\bf Lemma.}
\emph{ Assume $f \in JN_p^+(\Rn)$ and take $0<b<2^{-n}$. Then, given a cube $Q_0$, we have
\begin{align*}
|\{ x \in Q_0: &M_{Q_0}^{+,d}(f-f_{Q_0^{+,2}})^+(x) >  \lambda \}| 
\\ &\leq \frac{aK_f^+}{\lambda} 
|\{ x \in Q_0: M_{Q_0}^{+,d}(f-f_{Q_0^{+,2}})^+(x) > b\lambda\}|^{1/q},
\end{align*}
whenever $$b\lambda \geq \intav_{Q_0^+} (f-f_{Q_0^{+,2}})^+ \D x.$$
Here $a=4(1-2^nb)^{-1}$ and $q$ is the conjugate exponent of $p$.} 
\begin{proof}
Setting
$$
E_Q(\lambda):=\{x \in Q: M_{Q}^{+,d}(f-f_{Q_0^{+,2}})^+(x)>\lambda\},
$$
we may write the statement as
\begin{equation}\label{p4}
|E_{Q_0}(\lambda)| \leq \frac{4K_f^+}{(1-2^nb)\lambda} |E_{Q_0}(b\lambda)|^{1/q}.
\end{equation}
Consider the function $(f-f_{Q_0^{+,2}})^+$ and form the decomposition as above at level 
$b\lambda$ to obtain a family of pairwise non-overlapping dyadic subcubes with
\begin{equation*}%\label{p5}
E_{Q_0}(b\lambda)=\bigcup_j Q_j.
\end{equation*}
Since $b\lambda < \lambda$, we have
$E_{Q_0}(\lambda) \subset E_{Q_0}
(b\lambda)$. It now follows %from (\ref{p5}) 
that
\begin{equation}\label{p6}
E_{Q_0}(\lambda)=\bigcup_j E_{Q_j}(\lambda).
\end{equation}
We claim that for every $j$,
\begin{equation}\label{p7}
|E_{Q_j}(\lambda)| \leq \frac{2}{(1-2^nb)\lambda} \int_{Q_j \cup Q_j^+} (f-f_{Q_j^{+,2}})^+ \D x.
\end{equation}
Consider the functions $g_j:=(f-f_{Q_j^{+,2}})^+.$
To prove (\ref{p7}) it suffices to show that
\begin{equation}\label{p8}
E_{Q_j}(\lambda) \subset \{ x \in Q_j: M_{Q_j}^{+,d}g_j(x)>(1-2^nb)\lambda \}.
\end{equation}
Indeed, \eqref{p7} then follows at once from the weak type estimate \eqref{p3} applied to the functions $g_j$ with 
$\lambda$ replaced by $(1-2^nb)\lambda$.
Let $x \in E_{Q_j}(\lambda)$ for some $j$. Then there exists a 
dyadic subcube $Q$ of $Q_j$ containing $x$ and satisfying
$$\intav_{Q+} (f-f_{Q_0^{+,2}})^+ > \lambda$$
From \eqref{p2} we have 
$$
\intav_{Q_j^{+,2}} (f-f_{Q_0^{+,2}})^+ 
\leq 2^nb\lambda.
$$ 
Combining these, we obtain 
\begin{align*}
(1-2^nb)\lambda & <  \intav_{Q^+} (f-f_{Q_0^{+,2}})^+ \D x-\intav_{Q_j^{+,2}} (f-f_{Q_0^{+,2}})^+ \D x \\
& \leq \intav_{Q^+} (f-f_{Q_0^{+,2}})^+ \D x- \left(\intav_{Q_j^{+,2}} f-f_{Q_0^{+,2}}\D x \right)^+ \\
& = \intav_{Q^+} (f-f_{Q_0^{+,2}})^+-(f_{Q_j^{+,2}}-f_{Q_0^{+,2}})^+ \D x \\
& \leq \intav_{Q^+} (f-f_{Q_j^{+,2}})^+ \D x \\
&\leq M^{+,d}_{Q_j}g_j(x).
\end{align*}
Having now seen that (\ref{p7}) holds, 
we use \eqref{p6} and sum over all $j$ to obtain 
\begin{align*}
|E_{Q_0}(\lambda)| & = \sum_j |E_{Q_j}(\lambda)| \\
& \leq \frac{2}{(1-2^nb)\lambda} \sum_j \int_{Q_j \cup Q_j^+} (f-f_{Q_j^{+,2}})^+ \D x\\
& = \frac{2}{(1-2^nb)\lambda} \sum_j |Q_j|^{1/q}|Q_j|^{-1/q}  \int_{Q_j \cup Q_j^+} (f-f_{Q_j^{+,2}})^+ \D x\\
& \leq \frac{4K_f^+}{(1-2^nb)\lambda} \left( \sum_j |Q_j| \right)^{1/q},
\end{align*}
where the last inequality follows from the H\"older inequality 
and the definition of $K_f^+$.
Remembering also that $E_Q(b\lambda)=\bigcup_j Q_j$, 
we obtain the desired estimate.
\end{proof}

We now complete the proof of the theorem by iterating the previous lemma. 
Except for a few details, this a just a repetition of the argument used in 
\cite{AaBeKaYu}.

{\bf \emph{Proof of the Theorem.}} 
%We may assume $K_f^+=1$ and $f_{Q_0^{+,2}}=0$.
Using the same notation as in the proof of the lemma, we shall show
\begin{equation}\label{p9}
|E_{Q_0}(\lambda)| \leq C\left(\frac{K_f^+}{\lambda}\right)^p.
\end{equation}
Let us choose $$\lambda_0 :=\frac{2K_f^+}{b|Q_0|^{1/p}}$$ 
and assume $\lambda > \lambda_0$. 
Then take $N \in \Z_+$ such that 
\begin{equation}\label{p10}
b^{-N}\lambda_0 \leq \lambda < b^{-(N+1)}\lambda_0=
\frac{2b^{-(N+2)}K_f^+}{|Q_0|^{1/p}}.
\end{equation}
%The assumptions $K_f^+=1$ and $f_{Q_0^{+,2}}=0$ imply
By the definition of $K_f^+$, we have
\begin{equation}\label{p11}
%\frac{1}{b} 
\frac{1}{|Q_0|}\int_{Q_0 \cup Q_0^+} 
(f-f_{Q_0^{+,2}})^+ \D x \leq \frac{2K_f^+}{|Q_0|^{1/p}}=b\lambda_0. 
\end{equation}
In particular, this implies
$$
\frac{1}{b} \intav_{Q_0^+} 
(f-f_{Q_0^{+,2}})^+ \D x \leq 
\lambda_0 \leq b^{-1}\lambda_0 \leq \ldots \leq b^{-N} \lambda_0,
$$
allowing us to apply the previous lemma successively $N$ times to estimate
the left-hand side of (\ref{p9}) as follows:
\begin{align*}
&|E_{Q_0}(\lambda)| \\ &\leq |E_{Q_0}(b^{-N}\lambda_0)|\\
&\leq \frac{aK_f^+}{b^{-N}\lambda_0}\cdot\left(\frac{aK_f^+}{b^{-N+1}
\lambda_0}\right)^{1/q}\cdot
\ldots \cdot
\left(\frac{aK_f^+}{b^{-1}\lambda_0}\right)^{1/q^{N-1}}|E_{Q_0}(\lambda_0)|^{1/q^N} \\
& \leq \frac{aK_f^+}{b\lambda}\cdot\left(\frac{aK_f^+}{b^{2}\lambda}\right)^{1/q}\cdot \ldots \cdot 
\left(\frac{aK_f^+}{b^{N}\lambda}\right)^{1/q^{N-1}}\cdot\left(\frac{2}
{\lambda_0}{\int_{Q_0 \cup Q_0^+} (f-f_{Q_0^{+,2}})^+ \D x}\right)^{1/q^N},
\end{align*}
where the last inequality follows from the weak type estimate (\ref{p3})
and the first inequality in (\ref{p10}).
By the choice of $\lambda_0$ and (\ref{p11}) we further estimate
\begin{align*}|E_{Q_0}(\lambda)| &\leq 
\left(\frac{aK_f^+}{\lambda}\right)^{1+q^{-1}+\ldots+q^{-(N-1)}} \cdot
b^{-(1+2q^{-1}+\ldots+Nq^{-(N-1)})}\cdot
(2b |Q_0|)^{1/q^N}\\
&= \left(\frac{aK_f^+}{\lambda}\right)^{p-p/q^{N}}\cdot b^{-(1+2q^{-1}+
\ldots+Nq^{-(N-1)})+q^{-N}}\cdot
2^{1/q^N} \cdot |Q_0|^{1/q^N}.
\end{align*}
Since both $1+2q^{-1}+\ldots+Nq^{-(N-1)}$ and $p-p/q^{N}$
remain bounded as $N \to \infty$, we have
$$
|E_{Q_0}(\lambda)| \leq C|Q_0|^{1/q^{N}}\left(
\frac{1}{\lambda}\right)^{p-p/q^{N}}.
$$
Finally, we notice that from the second inequality in (\ref{p10}) 
we get
$$
|Q_0|^{1/q^{N}}\left(
\frac{1}{\lambda}\right)^{-p/q^{N}}=\lambda^{p/q^{N}}
|Q_0|^{1/q^{N}}\leq 2^{p/q^N}
b^{-(N+2)p/q^{N}}\leq C
$$
with $C$ independent of $N$.
Thus we have arrived at the desired estimate. 

For $0<\lambda\leq\lambda_0$ 
we use the trivial estimate
$$
|\{x \in Q_0: (f(x)-f_{Q_0^{+,2}})^+ > \lambda \}| \leq |Q_0| =\frac{2^p(K_f^+)^p}{b^p\lambda_0^p}
\leq C\left(\frac{K_f^+}{\lambda}\right)^p.
$$ 
\hfill{$\square$}

{\bf \emph{Acknowledgement.}} The author was supported 
by the Finnish Cultural Foundation. The author wishes to thank J. Kinnunen
 for proposing the problem.

\vspace{0.5cm}
\noindent
\small{\textsc{L. B.},}
\small{\textsc{Department of Mathematics},}
\small{\textsc{FI-90014 University of Oulu},}
\small{\textsc{Finland}}\\
\footnotesize{\texttt{lauri.berkovits@oulu.fi}}


\begin{thebibliography}{10}

\bibitem{AaBeKaYu}
D. Aalto, L. Berkovits, O. E. Kansanen and H. Yue
\newblock John-Nirenberg lemmas for a doubling measure. 
\newblock {\em Studia Math.} {\bf 204} (1) (2011), 21-37. 

%\bibitem{AiFoMa}
%H. Aimar, L. Forzani and F. J. Mart\'in-Reyes. 
%\newblock On weighted inequalities for 
%singular, integrals.
%\newblock {\em Proc. Amer. Math. Soc.} {\bf 125} (7) (1997), 2057-2064.

\bibitem{AiCr} 
H. Aimar and R. Crescimbeni. 
\newblock On one-sided $\bmo$ and Lipschitz functions. 
\newblock {\em Ann. Sc. Norm. Sup. Pisa} {\bf 27} (3-4) (1998), 437-456.

\bibitem{BaJiMa}
N. Badr, A. Jim\'enez-del-Toro and J. M. Martell.
\newblock $L^p$ self-improvement of generalized Poincar\'e inequalities in spaces of homogeneous type.
\newblock {\em J. Funct. Anal.} {\bf 260} (11) (2011), 3147-3188.


\bibitem{Be} 
L. Berkovits. 
\newblock Parabolic Muckenhoupt weights in the Euclidean space. 
\newblock {\em J. Math. Anal. Appl} {\bf 379} (2011), 524-537.


\bibitem{FrPeWh}
B. Franchi, C. P{\'e}rez and R. L. Wheeden,
\newblock Self-improving properties of John-Nirenberg and Poincar{\'e} inequalities on spaces of homogeneous type.
\newblock {\em J. Funct. Anal.}, 153 (1): 108--146, 1998.

\bibitem{GiMa}
M. Giaquinta and L. Martinazzi,
\newblock{\em An introduction to the regularity theory for elliptic systems, 
harmonic maps and minimal graphs.}
\newblock Edizioni della Normale, Pisa, 2005. 

\bibitem{Gi}
E. Giusti,
\newblock {\em Direct methods in the calculus of variations}.
\newblock World Scientific Publishing Co. Inc., River Edge, NJ, 2003.

\bibitem{Heik1}
T. Heikkinen,
\newblock Self-improving properties of generalized Orlicz-Poincar{\'e} inequalities.
\newblock {\em Rep. Univ. Jyväskylä Dept. Math. Stat.} 105, 2006.

\bibitem{JiMa}
A. Jim\'enez-del-Toro and J. M. Martell.
\newblock $L^p$ self-improvement of Poincar\'e type inequalities associated with approximations 
of the identity and semigroups.
\newblock Preprint.



\bibitem{JoNi}
F. John and L. Nirenberg.
\newblock On functions of bounded mean oscillation.
\newblock {\em Comm. Pure Appl. Math.} {\bf 14} (1961), 415-426.

%\bibitem{LeOm}
%A. Lerner and S. Ombrosi.
%\newblock A boundedness criterion for general maximal operators.
%\newblock {\em Publ. Mat.} {\bf 54} (1) (2010), 53-71.

\bibitem{LePe}
A. Lerner and C. P\'erez.
\newblock Self-improving properties of generalized Poincar\'e type 
inequalities throught rearrangements.
\newblock {\em Math. Scand.} {\bf 97} (2) (2005), 217-234.




\bibitem{MaPe}
P. MacManus and C. P{\'e}rez,
\newblock Generalized Poincar{\'e} inequalities: sharp self-improving properties.
\newblock {\em Internat. Math. Res. Notices}, 2:101--116, 1998.


%\bibitem{Ma1} 
%F. J. Mart\'in-Reyes
%\newblock New proofs of weighted inequalities for the one-sided 
%Hardy-Littlewood maximal functions.
%\newblock {\em Proc. Amer. Math. Soc.} {\bf 117} (1993), 691-698.

%\bibitem{MaPiTo} 
%F. J. Mart\'in-Reyes, L. Pick and A. De La Torre. 
%\newblock $A_\infty^+$ condition. 
%\newblock {\em Can. J. Math.} {\bf 45} (6) (1993) 1231-1244.


\bibitem{MaTo1} 
F. J. Mart\'in-Reyes and A. de la Torre.
\newblock One-sided $\bmo$ spaces. 
\newblock {\em J. London Math. Soc.} {\bf 49} (2) (1994), 529-542.


%\bibitem{MaTo2} 
%F. J. Mart\'in-Reyes and A. de la Torre.
%\newblock Two weight norm inequalities for one-sided fractional maximal operators. 
%\newblock {\em Proc. Amer. Math. Soc.} {\bf 117} (2) (1992), 483-489.


%\bibitem{MaOrTo} 
%F. J. Mart\'in-Reyes, P. Ortega Salvador and A. de la Torre.
%\newblock Weighted 
%inequalities for one-sided maximal functions.
%\newblock {\em Trans. Amer. Math. Soc.} {\bf 319} (2) (1990), 517-534.

%\bibitem{Mo1} 
%J. Moser. 
%\newblock A Harnack inequality for parabolic differential equations.
%\newblock {\em Comm. Pure Appl. Math.} {\bf 17} (1964), 101-134.

%\bibitem{Mo2} 
%J. Moser.
%\newblock Correction to: ``A Harnack inequality for parabolic differential equations''.
%\newblock {\em Comm. Pure Appl. Math.} {\bf 20} (1967), 231-236.

%\bibitem{Om1} 
%S. Ombrosi.
%\newblock Weak weighted inequalities for a dyadic one-sided maximal function 
%in $\Rn$.
%\newblock {\em Proc. Amer. Math. Soc.} {\bf 133} (2005), 1769-1775.

%\bibitem{Sa1} 
%E. Sawyer.
%\newblock Weighted inequalities for the one-sided Hardy-Littlewood maximal 
%functions.
%\newblock {\em Trans. Amer. Math. Soc.} {\bf 297} (1) (1986), 53-61.


%\bibitem{St1}
%E. M. Stein.
%\newblock {\em Harmonic analysis: real-variable methods, orthogonality, and
% oscillatory integrals}, volume~43 of {\em Princeton Mathematical Series}.
%\newblock Princeton University Press, Princeton, NJ, 1993.
%\newblock With the assistance of Timothy S. Murphy, Monographs in Harmonic Analysis, III.

\bibitem{To}
A. Torchinsky,
\newblock {\em Real-variable methods in harmonic analysis}.
  Pure and Applied Mathematics, 123.
\newblock Academic Press Inc., Orlando, FL, 1986.
\end{thebibliography}
\end{document}